\documentclass[12pt,twoside]{article}
\usepackage{hyperref}
\usepackage{amsthm,amssymb,epsfig,amscd,amsmath,color,graphicx}
\newtheorem{theorem}{Theorem}[section]
\newtheorem{property}{Property}
\newtheorem*{MainTheorem}{Main Theorem}
\newtheorem*{MainLemma}{Main Lemma}

\newtheorem{example}[theorem]{Example}

\newtheorem{open problem}[theorem]{Open Problem}

\newtheorem{definition}[theorem]{Definition}

\newtheorem{lemma}[theorem]{Lemma}

\newcommand{\inv}{^{-1}}

\newcommand{\R}{\mathbb R}

\title{Isometries of Products of Path-Connected Locally Uniquely Geodesic Metric Spaces with the Sup Metric are Reducible.}
\author{William Malone}

\begin{document}

\maketitle

\begin{abstract}
Let $M_i$ and $N_i$ be path-connected locally uniquely geodesic metric spaces that are not points and $f:\prod_{i=1}^m M_i\rightarrow \prod_{i=1}^n N_i$ be an isometry where $\prod_{i=1}^n N_i$ and $\prod_{i=1}^m M_i$ are given the sup metric.  Then $m=n$ and after reindexing $M_i$ is isometric to $N_i$ for all $i$.  Moreover $f$ is a composition of an isometry that reindexes the factor spaces and an isometry that is a product of isometries $f_i:M_i\rightarrow N_i$.  \end{abstract}

\section{Introduction}

For the duration a product $\prod_{i=1}^m M_i$ of metric spaces is always considered to be endowed with the sup metric.  In these metric spaces there are two obvious types of isometries $f:\prod_{i=1}^m M_i\rightarrow \prod_{i=1}^m M_i$.  The first is a product $f=f_1\times\dots\times f_m$ where $f_j:M_j\rightarrow M_j$ is an isometry.  The second is a reindexing of the form $$g(\alpha_1,\dots,\alpha_m)=(\alpha_{\pi\inv(1)},\dots,\alpha_{\pi\inv(m)})$$ for some permutation $\pi$ of $\{1,\dots,m\}$.  This leads to the following definition.

\begin{definition}
 An isometry $f:\prod_{i=1}^m M_i\rightarrow \prod_{i=1}^m N_i$ is reducible if there is a permutation $\pi$ of $\{1,2,\dots,m\}$ and isometries $f_i:M_i\rightarrow N_{\pi(i)}$ such that $$f(\alpha_1,\dots,\alpha_m)=(f_{\pi\inv(1)}(\alpha_{\pi\inv(1)}),\dots,f_{\pi\inv(m)}(\alpha_{\pi\inv(m)})).$$
\end{definition}

The main result states that in certain cases all isometries are reducible.

\begin{MainTheorem}
 Suppose $M_i$ and $N_i$ are path-connected locally uniquely geodesic metric spaces that are not points with $f:\prod_{i=1}^m M_i\rightarrow \prod_{i=1}^n N_i$ an isometry. Then $m=n$ and $f$ is reducible.\label{main}
\end{MainTheorem}

The heart of the proof lies in the fact that particular isometric embeddings of a graph in the product space $\prod_{i=1}^m M_i$ are invariant under any isometry.  This will quickly lead to the number of factors in the product being an isometry invariant as well as forcing $f$ to take the aspect of a product map.  A proof for the case where $m=2$ occurs in the appendix of \cite{FarbMasur08}, and a discussion of the background on geodesic metric spaces can be found in \cite{BriHae99}.

\subsection{Acknowledgements}  I would like to thank Mladen Bestvina, Yael Algom-Kfir, and Erika Meucci for their many helpful comments and discussions.  I would also like to thank Benson Farb for mentioning this problem during a talk at the University of Utah. 

\section{Proof of Main Theorem}

In this section we define a useful family of subsets of $\prod_{i=1}^m M_i$.  Utilizing the Main Lemma, a proof of the Main Theorem is then provided.

\begin{definition}
A $k$-slice $\Omega$ is a subset of $\prod_{i=1}^m M_i$ of the form $$\alpha_1\times  \dots \times \alpha_{k-1}\times \Omega_k\times \alpha_{k+1}\times \dots\times\alpha_m $$ where $\alpha_i\in M_i$ and $\Omega_k\subset M_k$ with $|\Omega_k|\geq 2$.  The $\alpha_i$ are called the fixed coordinates of $\Omega$.\label{slice}
\end{definition}

\begin{MainLemma} Suppose that $M_i$ and $N_i$ are path-connected locally uniquely geodesic metric spaces 
and $f:\prod_{i=1}^m M_i\rightarrow \prod_{i=1}^n N_i$ an isometry.  If $\Omega\subset \prod_{i=1}^m M_i$ 
is a $k$-slice then $f(\Omega)$ is a $j$-slice.  Moreover $\Omega$ and $\Omega'$ are both
 $k$-slices if and only if $f(\Omega)$ and $f(\Omega')$ are both $j$-slices. 
\end{MainLemma}

In order to prove the Main Lemma and hence the Main Theorem the fact that $m=n$ for any isometry $f:\prod_{i=1}^m M_i\rightarrow \prod_{i=1}^n N_i$ is essential.  However the proof of this fact is more technical and will be postponed.

\begin{proof}(of Main Theorem)
 Assuming $m=n$, let $\Phi\subset\prod_{i=1}^m M_i$ be a $k$-slice for $k\in\{1,\dots,m\}$.  
 By the Main Lemma the image $f(\Phi)$ is a $j$-slice for some $j$.  The second part of the Main Lemma gives more information namely that the image of any $k$-slice $\Omega\subset\prod_{i=1}^m M_i$ will be a $j$-slice.

To complete the proof it  suffices to show that if two points,  $$a=(\alpha_1,\alpha_2,\dots,\alpha_{k-1},\gamma,\alpha_{k+1},\dots,\alpha_m)$$ and $$b=(\beta_1,\beta_2,\dots,\beta_{k-1},\gamma,\beta_{k+1},\dots,\beta_m)$$
such that $\alpha_i\neq \beta_i$ for all $i$, have the same $k^{th}$ coordinate then $f(b)$ and $f(a)$ have the same $j^{th}$ coordinate.  This will show the existence of isometries $f_i:M_k\rightarrow N_j$ and complete the proof of reducibility after post composing by the reindexing isometry $\pi:\prod_{i=1}^mN_i\rightarrow \prod_{i=1}^m N_i$ where $\pi(N_j)=N_k$.  Consider the points $$\begin{array}{l} a_0=a=(\alpha_1,\alpha_2,\dots,\alpha_{k-1},\gamma,\alpha_{k+1},\dots,\alpha_m)\\ a_1=(\beta_1,\alpha_2,\dots,\alpha_{k-1},\gamma,\alpha_{k+1},\dots,\alpha_m)\\ a_2=(\beta_1,\beta_2,\alpha_3,\dots,\alpha_{k-1},\gamma,\alpha_{k+1},\dots,\alpha_m)\\
\vdots\\
a_{k-1}=(\beta_1,\beta_2,\dots,\beta_{k-1},\gamma,\alpha_{k+1},\dots,\alpha_m)\\
a_{k}=(\beta_1,\beta_2,\dots,\beta_{k-1},\gamma,\beta_{k+1},\alpha_{k+2},\dots,\alpha_m)\\
\vdots\\
a_{m-2}=(\beta_1,\beta_2,\dots,\beta_{k-1},\gamma,\beta_{k+1},\dots,\beta_{m-1},\alpha_m)\\
a_{m-1}=b=(\beta_1,\beta_2,\dots,\beta_{k-1},\gamma,\beta_{k+1},\dots,\beta_m).\\ \end{array}$$

Notice that $\{a_i,a_{i+1}\}$ is an $i+1$-slice for $i\in\{0,1,\dots,k-2\}$ and an $i+2$ slice for $i\in \{k-1,\dots,m-2\}$.  Since no pair $\{a_i,a_{i+1}\}$ is a $k$-slice, the pair of points $\{f(a_i),f(a_{i+1})\}$ is never a $j$-slice by the Main Lemma.  Thus the $j^{th}$ coordinate of $f(a_i)$ and $f(a_{i+1})$ are identical for all $i$.
\end{proof}

\section{The Graph $\mathcal Q^m_r$}

In this section we define a metric graph, give examples of particularly nice isometric embeddings of this graph in $\prod_{i=1}^m M_i$, show that $m=n$ for any isometry $f:\prod_{i=1}^m M_i\rightarrow \prod_{i=1}^n N_i$, and prove the Main Lemma.  

Observe that if $\alpha_i,\beta_i\in M_i$ are connected by a unique geodesic of length $r$ for $i=1,\dots,m$ then $a=(\alpha_1,\dots,\alpha_m)$ and $b=(\beta_1,\dots,\beta_m)$ are connected by a unique geodesic in $\prod M_i$.  This fact that in certain directions within the product pairs of  points are joined by unique geodesics is the main motivation for the following graph.  


\begin{definition}
 A quadrilateral graph of dimension $m$ with length $r$ denoted $\mathcal Q^m_r\subset \prod_{i=1}^m \R$ is the metric graph whose vertices are of the form $\pm e_i^r=(0,0,\dots,0,\pm 2r,0,\dots,0)$ (non-zero in $i^{th}$ coordinate) or $(\pm r,\pm r,\dots,\pm r)$.  Connect two vertices by an edge of length $r$ if their distance is $r$ in $\prod_{i=1}^m \R$.  
\end{definition}  

The first two examples of this family of graphs (pictured in Figure 1) are a quadrilateral with subdivided sides and the one skeleton of Kepler's Rhombic dodecahedron \cite{BallCox87}.

\begin{figure}[htb]
\begin{center}
\includegraphics[angle=270]{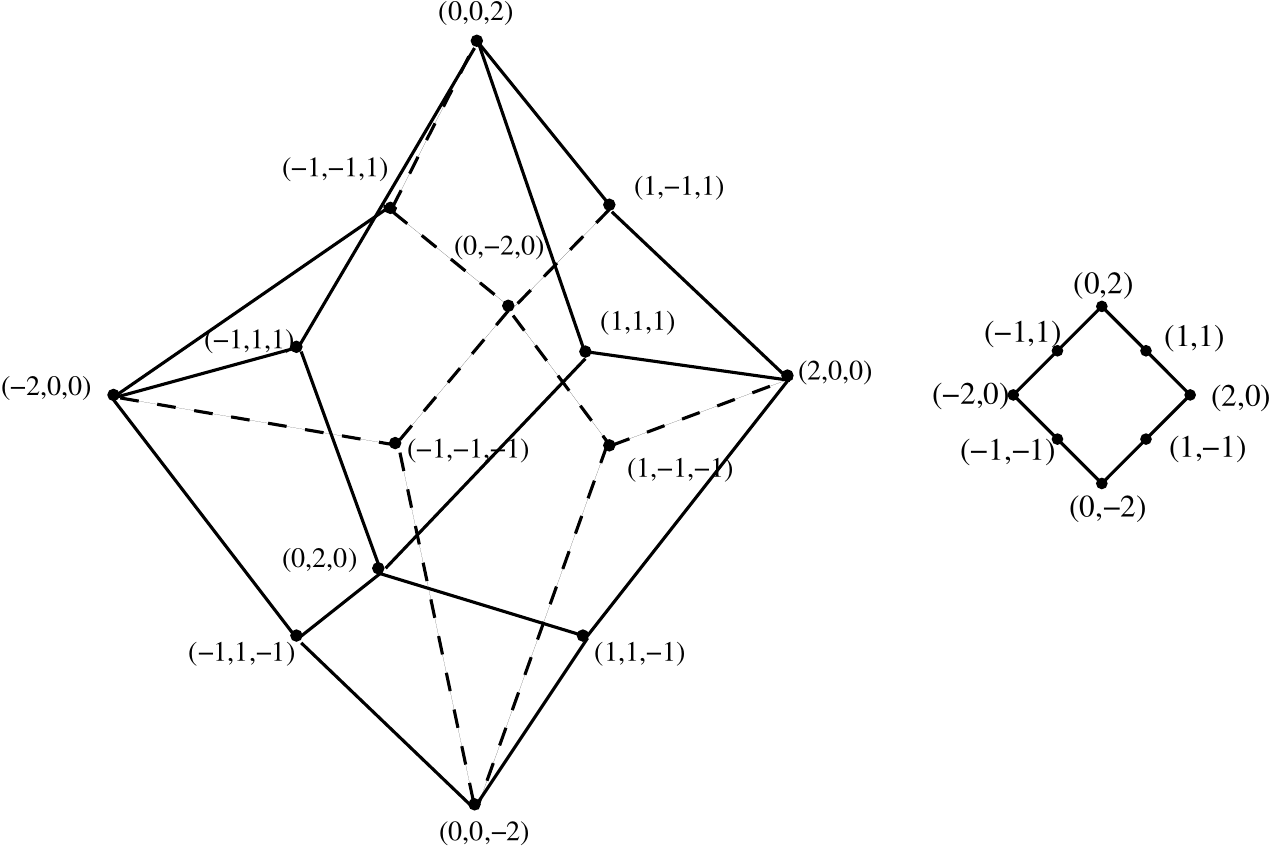}
 \caption{$\mathcal Q^3_1\subset \prod_{i=1}^3 \R$ (left) and $\mathcal Q^2_1\subset \prod_{i=1}^2 \R$ (right) }
\end{center}\label{graph}
\end{figure}

\begin{example}
 Let $M_i$ be a uniquely geodesic metric space.  An example of an isometric embedding $\iota:\mathcal Q^m_r\rightarrow \prod_{i=1}^m M_i$ is given by the following construction.  Pick a geodesic segment $\lambda_i\subset M_i$ of length $4r$.  Let $\lambda_i(0)=\alpha_i$, $\lambda_i(r)=\theta_i$ $\lambda_i(2r)=\beta_i$, $\lambda_i(3r)=\varphi_i$ and $\lambda_i(4r)=\omega_i$.  Then let $$\begin{array}{ll}& \iota(e^r_k)=\{(\beta_1,\beta_2,\dots,\beta_{k-1},\omega_k,\beta_{k+1},\dots,\beta_m\|k\in \{1,2,\dots,m\}\}\\
&\iota(-e^r_k)=\{(\beta_1,\beta_2,\dots,\beta_{k-1},\alpha_k,\beta_{k+1},\dots,\beta_m\|k\in \{1,2,\dots,m\}\}\\                                                                                                        
 and\,\,as \,\,a \,\,set& \iota(\pm r, \pm r,\dots,\pm r)=\{(\pi_1,\dots,\pi_m)|\pi_i\in\{\theta_i,\varphi_i\}\}
\end{array}$$  
which defines the isometry on the set of vertices.  One can check that appropriate vertices are connected by unique geodesics of length $r$.  \label{graph}
\end{example}

This construction shows that isometric embeddings $\iota:\mathcal Q^m_r\rightarrow \prod_{i=1}^m M_i$ are easy to construct when the $M_i$ are uniquely geodesic metric spaces.  If the spaces $M_i$ are only path-connected locally uniquely geodesic metric spaces then in general an isometric embedding $\iota:\mathcal Q^m_r\rightarrow \prod_{i=1}^m M_i$ will not exist for large $r$.  However for every $a\in \prod M_i$ there exists a metric ball $B_{\epsilon}(a)$ such that $B_{\epsilon}(a)$ is a product of uniquely geodesic metric spaces.  Thus for $4r<\epsilon$ an isometrically embedded $\mathcal Q^m_r$ can be constructed in the same way.

\begin{definition}
Let $\iota:\mathcal Q^k_r\rightarrow \prod_{i=1}^m M_i$ be an isometric embedding.  If the image of every edge is a uniquely geodesic segment in $\prod M_i$ then call $\iota:\mathcal Q^k_r\rightarrow \prod_{i=1}^m M_i$ admissible.
\end{definition}

If $\iota:\mathcal Q^k_r\rightarrow \prod_{i=1}^m M_i$ is admissible and $f:\prod_{i=1}^m M_i\rightarrow \prod_{i=1}^n N_i$ is an isometry then $f\circ \iota:\mathcal Q^k_r \rightarrow \prod_{i=1}^n N_i$ is admissible.  This is clear from definitions. 

\begin{definition}
 Let $\iota:\mathcal Q^k_r\rightarrow \prod_{i=1}^m M_i$ be admissible.  The map $\iota$ is standard if for all $j$ there exists an $l$ such that $\{\iota(e^r_j),\iota(-e^r_j)\}$ is an $l-$slice. \end{definition}

The proof of the Main Lemma requires a result establishing the fact that every admissible isometric embedding $\iota:\mathcal Q^m_r\rightarrow \prod_{i=1}^m M_i$ is standard.  The following discussion works toward establishing this fact.  First some results concerning the behavior of points in the image of  an admissible $\iota:\mathcal Q^m_r\rightarrow \prod_{i=1}^m M_i$ are required.

\begin{property}
For any vertex $x_j\in Q^m_r$ of the form $\{\pm e^r_1,\dots,\pm e^r_m\}$ and any pair of points $\{e_i^r,-e_i^r\}$ with $x_j\not \in \{e_i^r,-e_i^r\}$ there exists a geodesic segment $\Lambda\subset \mathcal Q^m_r$ with endpoints $\{e_i^r,-e_i^r\}$ such that $x_j\in \Lambda$. \label{prop1}
\end{property}

Constructing such a piecewise uniquely geodesic path in $\mathcal Q^m_r$ is straightforward.    Note that there are precisely $2^{2m-4}$ paths connecting pairs of endpoints $\{e_i^r,-e_i^r\}$ so $\Lambda$ is unique only when $m=2$.

\begin{property}  Assume that $M_i$ is a uniquely geodesic metric space.  Let $p$ and $q$ be points in $\prod_{i=1}^m M_i$ and $\Lambda$ a geodesic with endpoints $\{p,q\}$.  Let $d_{M_i}(p_i,q_i)$ be the distance between the $i^{th}$ coordinates of $p$ and $q$. Note that $d_{M_i}(p_i,q_i)=d(p,q)$ if and only if the $i^{th}$ coordinate path of $\Lambda$ (denoted $\lambda_i$) is a geodesic of length $d(p,q)$.  Moreover for any point $z\in \Lambda$ the $i^{th}$ coordinate of $z$ is uniquely determined by $d(p,z)$ (since $\lambda_i$ is the only geodesic connecting its endpoints in $M_i$).
\label{prop2} \end{property}

Property \ref{prop2} follows directly from the definition of a geodesic and will allow us to characterize the possible admissible $\iota:\mathcal Q^m_r\rightarrow \prod_{i=1}^m M_i$ when $M_i$ is uniquely geodesic.

\begin{lemma}
 Assume that $M_i$ is a uniquely geodesic metric space and let $\iota:\mathcal Q^k_r\rightarrow \prod_{i=1}^m M_i$ be admissible.  Let $$\begin{array}{c}\iota(e_j^r)=(\alpha_1,\dots,\alpha_m)\\ \iota(-e_j^r)=(\beta_1,\dots,\beta_m)\\ \iota(e_t^r)=(\gamma_1,\dots,\gamma_m)\\ \iota(-e_t^r)=(\kappa_1,\dots,\kappa_m)\\ \end{array}.$$  If $d(\iota(e_t^r),\iota(-e_t^r))=d_{M_l}(\gamma_l,\kappa_l)=4r$  then $\alpha_l=\beta_l$.  \label{disjoint}
\end{lemma}

\begin{proof} 
By property \ref{prop1} there exists geodesics $\Lambda,\Lambda'\subset \iota(\mathcal Q^k_r)$ such that the endpoints of $\Lambda$ and $\Lambda '$ are $\iota(\pm e_t^r)$ with $\iota(e_j^r)\in \Lambda$ and $\iota(-e_j^r)\subset \Lambda '$.  By property \ref{prop2} and the fact that $d(\iota(e_t^r), \iota(e_j^r))=d(\iota(e_t^r), \iota(-e_j^r))=2r$ we see that $\alpha_l=\beta_l $.
\end{proof}

For an admissible $\iota:\mathcal Q^k_r\rightarrow\prod_{i=1}^m M_i$ associate a number $$q_j=|\{l\in\{1,\dots,m\} | d(\iota(e_j^r),\iota(-e_j^r))=d_{M_l}(\alpha_l,\beta_l)\}|$$ to each pair $\iota(\pm e_j^r)$.  Clearly $q_j\geq 1$ for all $j$.  If follows immediately from Lemma \ref{disjoint} that two pairs $\{\iota(e_j^r),\iota(-e_j^r)\}$ and $\{\iota(e_t^r),\iota(-e_t^r)\}$ cannot attain the distance $4r$ in the same coordinate.  These two observations imply $q_j=1$ for all $j$ when $k=m$.  

\begin{lemma}
  Let $M_i$ and $N_i$ be path-connected locally uniquely geodesic metric spaces and $f:\prod_{i=1}^m M_i \rightarrow \prod_{i=1}^n N_i$ be an isometry.  Then $m=n$.\label{m=n}
\end{lemma}

\begin{proof}(Of lemma \ref{m=n})
 By invariance of an admissible $\iota(\mathcal Q^k_r)$ under isometry it follows that the largest $k$ such that there exists an admissible $\iota:\mathcal Q^k_r\rightarrow \prod_{i=1}^m M_i$ is an isometry invariant.  Let this number be denoted by $L(\prod_{i=1}^m M_i)$.  There exists a ball $B(a,4r)$ around any point $a\in \prod_{i=1}^m M_i$ such that $B(a,4r)$ is a product of uniquely geodesic metric spaces.  An admissible $\mathcal Q^m_r$ can be constructed in $B(a,4r)$ by Example 3.2 which shows $L(\prod_{i=1}^m M_i)\geq m$.  Since the number $q_j\geq1$ for all pairs $\{\iota(e_j^r),\iota(-e_j^r)\}$ in any admissible $\iota:\mathcal Q^k_r\rightarrow \prod_{i=1}^m M_i$ and $\sum_{j=1}^m q_j\leq m$ it is obvious that for $k>m$ no such $\iota$ can exist and thus $L(\prod_{i=1}^m M_i)\leq m$.  
\end{proof}

\begin{lemma}
 Let $M_i$ be a uniquely geodesic metric space and $\iota:\mathcal Q^m_r \rightarrow \prod_{i=1}^m M_i$ admissible.  Then  $\iota:\mathcal Q^m_r\rightarrow \prod_{i=1}^m M_i$ is standard.\label{standard}
\end{lemma}

\begin{proof} (Of Lemma \ref{standard})  Using the above fact that $q_j=1$ for all $j$ and Lemma \ref{disjoint}  the collection of pairs $$\left(\cup_{i=1}^{j-1} \{\iota(e_i^r),\iota(-e_i^r)\}\right)\cup\left(\cup_{i=j+1}^n\{\iota(e_i^r),\iota(-e_i^r)\}\right)$$ forces the distance between coordinates of $\{\iota(e_j^r),\iota(-e_j^r)\}$ to be either $4r$ or $0$,  thus making $\{\iota(e_j^r),\iota(-e_j^r)\}$ a $l$-slice for some $l$. \end{proof}

We now restate and prove the Main Lemma.

\begin{MainLemma} Suppose that $M_i$ and $N_i$ are path-connected locally uniquely geodesic metric spaces 
and $f:\prod_{i=1}^m M_i\rightarrow \prod_{i=1}^n N_i$ an isometry.  If $\Omega\subset \prod_{i=1}^m M_i$ 
is a $k$-slice then $f(\Omega)$ is a $j$-slice.  Moreover $\Omega$ and $\Omega'$ are both
 $k$-slices if and only if $f(\Omega)$ and $f(\Omega')$ are both $j$-slices. 
\end{MainLemma}

\begin{proof}
By lemma \ref{m=n} $m=n$.  Three things must be shown.
 
\begin{enumerate}
\item[Step 1] (Part A) Given three distinct points $a,b,c$ such that $\{a,b\}$ is a $j$-slice, $\{b,c\}$ is a $k$-slice, and $\{a,c\}$ is an $i$-slice then $i=j=k$.  

(Part B) Given 4 distinct points $a,b,c,d$ such that $\{a,b\}$, $\{b,c\}$, $\{c,d\}$, and $\{a,d\}$ are $i,j,k,$ and $l$-slices respectively then either $i=j=k=l$ or $i=k$ and $j=l$.

\item[Step 2] If $\Omega$ is a $k$-slice then $f(\Omega)$ is a $j$-slice.

\item[Step 3] If $\Omega$ and $\Omega'$ are $k$-slices then $f(\Omega)$ and $f(\Omega')$ are $j$-slices.

\end{enumerate}

\noindent(Proof of Step 1 Part A)  By definition of a $k$-slice if $i=j$ then $i=j=k$.  Suppose $i\neq j$ then $a$ and $c$ differ in the $i^{th}$ and $j^{th}$ coordinate.  But $\{a,c\}$ is a $k$-slice which is a contradiction.   \vspace{5mm}

\noindent(Proof of Step 1 Part B)  Observe that since $\{a,b\}$ is a $i$-slice and $\{b,c\}$ is a $j$-slice the points $a$ and $c$ differ in the $i^{th}$ and $j^{th}$ coordinates.  Similarly they differ in the $l^{th}$ and $k^{th}$ coordinate.  Thus $\{i,j\}=\{k,l\}$.  If $i=l$ then $i=j=k=l$ otherwise $i=k$ and $j=l$.\vspace{5mm}

\noindent(Proof of Step 2)  
Any $k$-slice $\Omega$ is contained in a $k-$slice of the form $$\Phi=\alpha_1\times\dots\alpha_{k-1} \times M_k\times \alpha_{k+1}\dots \times\alpha_m$$ for $\alpha_i\in M_i$ so it suffices to show the result for $\Phi$.   Let $\eta\in M_k$ be an arbitrary point.  Let $\delta_s$ and $\bar{\epsilon_{\eta}}$ be the uniquely geodesic constants for $\alpha_s\in M_s$ and $\eta\in M_k$ respectively.  Define $$W_{\eta}=\left(\prod_{s=1}^{k-1} B(\alpha_s,\epsilon_{\eta})\right)\times B(\eta_t,\epsilon_{\eta})\times \left(\prod_{s=k+1}^{m} 
   B(\alpha_s,\epsilon_{\eta})\right)$$ where $\epsilon_{\eta}=min\{\bar{\epsilon_{\eta}},\delta_1,\dots,\delta_m\}$.  Let $\Psi=W_{\eta}\cap \Phi$.  For any three points $x_1,x_2,x_3\in  \Psi$ construct three isometrically embedded graphs $\iota_{i}: \mathcal Q^m_r\rightarrow W_{\eta}$ with $r=\frac{1}{4}d(x_i,x_{i+1})$, $\iota_i(e_1^r)=x_i$, and $\iota_i(-e_1^r)=x_{i+1}$  for $i\in \{1,2,3\}$ via Example 3.2 (note that $x_4=x_1$).  Since $W_{\eta}$ is a product of uniquely geodesic metric spaces we can apply lemma \ref{standard} to $f\circ \iota_i$ to see that $\{f(x_1),f(x_2),f(x_3)\}$ satisfy the hypothesis of Step 1 part A.  Thus $\{f(x_1),f(x_2),f(x_3)\}$ is a $j$-slice and since $x_i$ were arbitrary $f(\Psi)$ is a $j$-slice.  If $\Psi'$ is a defined about a point $\eta'\in M_k$ sufficiently close to $\eta$ then $f(\Psi)\cap f(\Psi')$ consists of more than two points and the set $f(\Psi)\cup f(\Psi')$ is a $j$-slice.   This implies $f(\Phi)$ is a $j-slice$ since $M_k$ has a single path component.\vspace{5mm}

 \noindent(Proof of Step 3)  Note that $$\Omega\subset \alpha_1\times\dots\times\alpha_{k-1} \times M_k\times\alpha_{k+1}\times\dots\times\alpha_m=\Phi$$ and $$\Omega'\subset \beta_1\times\dots\times\beta_{k-1} \times M_k\times\beta_{k+1}\times\dots\times\beta_m=\Phi'$$ for some $M_k$ and
hence it suffices to show the result for $\Phi,\Phi'$.  We can also assume that $\alpha_i=\beta_i$ for all $i\neq 1$ since we can change one coordinate at a time to interpolate between $\Phi$ and $\Phi'$.  Now there clearly exist distinct points  $w_1,w_2\in \Phi$ and
$w_1',w_2'\in \Phi'$ such that $\{w_1,w_1'\}$ and $\{w_2,w_2'\}$ are $1$-slices.  By step 2 $f(w_1),f(w_1'),f(w_2),f(w_2')$ satisfy the hypothesis of Step 1 part B.  If $\{f(w_1),f(w_1'),f(w_2),f(w_2')\}$ was a $j$-slice then by step 2 after applying the inverse isometry $\{w_1,w_2,w_1,'w_2'\}$ would be a $k$-slice or a $1$-slice which is a contradiction.
\end{proof}

\bibliography{prodRM8}

\noindent William Malone:

\noindent Dept. of Mathematics, University of Utah

\noindent 155 S 1400 E Room 233,

\noindent Salt Lake City, UT 84112-0090

\noindent Email:  malone@math.utah.edu

\end{document}